\documentclass[conf]{new-aiaa}
\usepackage{hyperref}

\usepackage[utf8]{inputenc}

\usepackage{graphicx}
\usepackage{amsmath}
\usepackage[version=4]{mhchem}
\usepackage{siunitx}
\usepackage{longtable,tabularx}
\usepackage{graphicx}
\graphicspath{{Figures/}}
\usepackage{dcolumn}
\usepackage{bm}

\usepackage{float}
\usepackage{caption}
\usepackage{subcaption}

\usepackage{soul}
\usepackage{amsmath,amssymb,amsfonts}
\usepackage{graphics}
\usepackage{color}
\usepackage{xcolor}
\definecolor{rev}{rgb}{0,0,1}

\usepackage{enumitem}

\usepackage{tabularx}
\usepackage{stackengine}

\title{Reduced order modeling of fluid flows: Machine learning, Kolmogorov barrier, closure modeling, and partitioning}

\author{Shady E. Ahmed \footnote{PhD Candidate, School of Mechanical \& Aerospace Engineering, Oklahoma State University, Stillwater, OK 74078, USA.},
Suraj Pawar \footnote{PhD Candidate, School of Mechanical \& Aerospace Engineering, Oklahoma State University, Stillwater, OK 74078, USA.},
Omer San \footnote{Assistant Professor, School of Mechanical \& Aerospace Engineering, Oklahoma State University, Stillwater, OK 74078, USA.} }
\affil{School of Mechanical \& Aerospace Engineering, Oklahoma State University, Stillwater, OK 74078, USA.}

\author{Adil Rasheed\footnote{Professor, Department of Engineering Cybernetics, Norwegian University of Science and Technology, N-7465, Trondheim, Norway.}}
\affil{Department of Engineering Cybernetics, Norwegian University of Science and Technology, N-7465, Trondheim, Norway.}

\begin{document}

\maketitle

\begin{abstract}
In this paper, we put forth a long short-term memory (LSTM) nudging framework for the enhancement of reduced order models (ROMs) of fluid flows utilizing noisy measurements. We build on the fact that in a realistic application, there are uncertainties in initial conditions, boundary conditions, model parameters, and/or field measurements. Moreover, conventional nonlinear ROMs based on Galerkin projection (GROMs) suffer from imperfection and solution instabilities due to the modal truncation, especially for advection-dominated flows with slow decay in the Kolmogorov width. In the presented LSTM-Nudge approach, we fuse forecasts from a combination of imperfect GROM and uncertain state estimates, with sparse Eulerian sensor measurements to provide more reliable predictions in a dynamical data assimilation framework. We illustrate the idea with the viscous Burgers problem, as a benchmark test bed with quadratic nonlinearity and Laplacian dissipation. We investigate the effects of measurements noise and state estimate uncertainty on the performance of the LSTM-Nudge behavior. We also demonstrate that it can sufficiently handle different levels of temporal and spatial measurement sparsity. This first step in our assessment of the proposed model shows that the LSTM nudging could represent a viable realtime predictive tool in emerging digital twin systems.
\end{abstract}




\section{Introduction} \label{sec:intro}
Reduced order modeling (ROM) is a family of protocols that aim at representing the system's dynamics of interest with minimal computational burden \cite{bai2002krylov,lucia2004reduced,hess2019localized,kramer2019nonlinear,swischuk2019projection,bouvrie2017kernel,hamzi2019local,korda2018data,korda2018linear,hartmann2018model,holmes2012turbulence,taira2017modal,taira2019modal,noack2011reduced,rowley2017model,nair2019transported,kaiser2014cluster,haasdonk2011training}. Standard approaches usually consist of two major steps; (1) tailor a low-order subspace, where the flow trajectory can be sufficiently approximated to live, (2) build a surrogate model to cheaply evolve this trajectory in time. For the former, modal decomposition techniques have shown substantial success in extracting the physically important features and underlying patterns of the flow. Examples include proper orthogonal decomposition (POD) \cite{sirovich1987turbulence,berkooz1993proper,holmes2012turbulence,chatterjee2000introduction,rathinam2003new}, balanced proper orthogonal decomposition (BPOD) \cite{willcox2002balanced,rowley2005model,singler2009proper,singler2012balanced}, spectral proper orthogonal decomposition (SPOD) \cite{sieber2016spectral,picard2000pressure,taira2017modal,towne2018spectral,schmidt2020guide}, and dynamic mode decomposition (DMD) \cite{schmid2010dynamic,kutz2016dynamic,chen2012variants,askham2018variable,rowley2009spectral,tu2014dynamic,mezic2013analysis,schmid2011applications}. Of particular interest, POD has gained historical recognition in fluid dynamics community, representing a set of data with minimal number of basis functions or modes while preserving as much energy as possible \cite{lumley1967structure,berkooz1993proper,holmes2012turbulence}. In particular, POD generates a set of hierarchically arranged modes, sorted by their respective contribution to the total variance of information in the data. In fluid flow applications, with velocity field data, this information corresponds to the flow's kinetic energy.

As a mathematical measure of the system's reducibility and the quality of a constructed (linear) subspace, the Kolmogorov $n$-width \cite{kolmogoroff1936uber} is a classical concept from approximation theory that quantifies the worst-case scenario error that might arise from the projection of solution trajectory onto an optimal subspace. Mathematically, it is defined as follows \cite{taddei2020registration,greif2019decay,pinkus2012n},
\begin{equation}
    d_n({\cal{M}}) := \inf_{{\cal{S}}_n} \sup_{q \in {\cal{M}}} \inf_{w\in {\cal{S}}_n} \| q - w \| ,
\end{equation}
where ${\cal{S}}_n$ is a linear $n$-dimensional subspace, and ${\cal{M}}$ is the solution manifold. The first infimum is taken over all possible $n$-dimensional subspaces, $q$ is a state on the solution manifold ${\cal{M}}$, while the last infimum sweeps all corresponding states that live in ${\cal{S}}_n$. In other words, $d_n({\cal{M}})$ estimates the largest error that might arise from approximating the solution manifold using the best-possible $n$-dimensional linear subspace. Assuming an orthogonal projection of $q$ onto ${\cal{S}}_n$ is possible, then the previous relation reduces to 
\begin{equation}
    d_n({\cal{M}}) := \inf_{{\cal{S}}_n} \sup_{q \in {\cal{M}}} \| q- \Pi_{{\cal{S}}_n} q \| ,
\end{equation}
where $\Pi_{{\cal{S}}_n}$ is the orthogonal projector onto ${\cal{S}}_n$. Using information about the decay rate of $d_n({\cal{M}})$ with increasing $n$, the system's reducibility and the quality of a reduced order approximation can be judged. Unfortunately, most of the fluid systems of practical relevance exhibit a slow decay of the Kolmogorov width, hindering reasonable approximation of the flow dynamics using a linear subspace. This has been denoted as the ``\emph{Kolmogorov barrier}'' in ROM context. Recently, efforts have been devoted to break-up or bypass this barrier by building more representative and concise subspaces. This can achieved either by partitioning techniques with the aim of localizing the resulting basis functions \cite{ijzerman2000signal,borggaard2007interval,san2015principal,ahmed2018stabilized,ahmed2019memory,ahmed2020breaking,amsallem2012nonlinear,peherstorfer2015online,taddei2015reduced,xiao2019domain} and preventing modal deformation, or constructing nonlinear latent subspaces using auto-encoders \cite{lee2020model,maulik2020reduced,gonzalez2018deep,kashima2016nonlinear,wiewel2019latent,masci2011stacked}.

Building surrogate models to evolve on reduced manifolds has been traditionally categorized into two groups; physics-based and data-based. Physics-based models rely on the governing equations from first principles, where the full order model (FOM) operators are projected (e.g., using Galerkin-type techniques) onto the reduced subspace to structure a reduced order model (ROM). Those are favorable because of their interpretability and generalizability, as well as the existence of robust techniques for stability and uncertainty analysis. However, those are usually expensive to represent for turbulent and advection-dominated flows, with a slow decay of the Kolmogorov $n$-width, necessitating an increase in the number of modes (or degrees of freedom) to be retained in ROM. Otherwise, modal truncation results in a Galerkin-based ROM (GROM) that might be eventually unstable. In this regard, closure techniques have been introduced to stabilize GROMs and account for the effect of truncated modes on the retained modes' dynamics \cite{wang2011two,borggaard2011artificial,iliescu2014variational,xie2017approximate,xie2018numerical,rahman2019dynamic,sirisup2004spectral,san2014proper,san2014basis,protas2015optimal,cordier2013identification,osth2014need,couplet2003intermodal,kalb2007intrinsic,kalashnikova2010stability,xie2018data,mohebujjaman2019physically,wang2012proper,akhtar2012new,balajewicz2012stabilization,amsallem2012stabilization,san2015stabilized,gunzburger2019evolve}. On the other hand, data-based models solely depend on archival data (thus called nonintrusive) to learn the underlying relations that govern the dynamical evolution of the system. Nonintrusive ROMs have benefited from the widespread machine learning (ML) tools to build stable and accurate models, compared to their GROM counterparts. In particular, (deep) neural networks have been extensively utilized to emulate the dynamical evolution of ROMs \cite{kutz2017deep, brunton2019machine, brenner2019perspective, duraisamy2019turbulence, xie2019non, jian2019flowfield,pawar2019deep,san2019artificial,rahman2019nonintrusive,maulik2020time,renganathan2020machine}. However, those often lack human interpretability and generalizability, and can even become prohibitively ``data-hungry''.

More recently, there has been a momentum in the research community to establish hybrid frameworks, that exploit machine learning algorithms and abundant data streams along with physical models to maximize their benefits \cite{reichstein2019deep, karpatne2017theory,rahman2018hybrid, san2018neural, wan2018data, xie2018data, mohebujjaman2019physically,maulik2019accelerating,discacciati2020controlling,pawar2020data,pawar2020evolve}. It has been shown that such hybridization can provide models that are superior to their individual components. Likewise, physics-informed machine learning is also gaining tremendous popularity, using human-knowledge and physical intuition to constrain the neural network predictions \cite{raissi2019physics,jagtap2020conservative,mao2020physics,yang2020physics}. Along similar lines, in this paper, we propose a hybrid framework that blends live measurement streams with physical models in order to achieve better predictions. Moreover, we suppose that both components (i.e., the physics-based model and data) are imperfect, thus avoiding biases in predictions. The physics-based ROMs (e.g., GROMs) are inherently imperfect due to the modal truncation and intrinsic nonlinearity. We also perturb the initial conditions to further mimic erroneous state estimates in practice. Meanwhile, we realize that, more often than not, sensor signals are noisy. Thus, we utilize recurrent neural networks, namely the long short-term memory (LSTM) variant, to combine the possibly defective model prediction with noisy measurements to ``nudge''  the model's solution towards the true states. 

Nudging is a data assimilation (DA) technique, which is a well-established family of predictive tools in geosciences, especially numerical weather forecasts \cite{ghil1991data,kalnay2003atmospheric,lewis2006dynamic,lorenc1991meteorological, derber1989global}. It works by relaxing the model state toward observations by adding correction (or nudging) terms, proportional to the difference between observations and model state, known as innovation in DA context (e.g., \cite{anthes1974data,lei2015nudging}). Usually, this proportionality is assumed to be linear, and the proportionality constants (or weights) are empirically tuned. Here, we use a simplistic LSTM architecture to generalize this relation to consider nonlinear mappings among the innovation and nudging terms. We apply the proposed framework (called LSTM-Nudge in the paper) for the reduced order modeling of the one-dimensional viscous Burgers equation as a starting benchmark problem for advection-dominated fluid flows with quadratic nonlinearity. We test the performance of LSTM-Nudge with various levels of measurement noises, initial field perturbations, and sensors signals sparsity. Therefore, the hybrid modeling approach presented in this paper illustrates a novel way of combining the best of both the physics-driven and the data-driven modeling approaches. Despite being shown in the context of a relatively simple dynamical system, this approach would be ideal for accurate and realtime modeling of complex systems, and therefore can be considered as a viable enabler for emerging digital twin technologies \cite{rasheed2020digital,ganguli2020digital,tao2018digital}. 

\section{Reduced Order Modeling} \label{sec:rom}
In this section, we present the reduced order formulations adopted in this study. In particular, we utilize proper orthogonal decomposition (POD) as a data-driven tool to extract the flow's coherent structures and build a reduced order subspace that best approximates the flow fields of interest. Then, we adopt a Galerkin approach to project the full order model operators onto that reduced space to build a ``'physics-constrained'' reduced order model.

\subsection{Governing equation} \label{sec:gov}
Here, we consider the one-dimensional (1D) viscous Burgers equation as a test bed. It represents a simple form of the Navier-Stokes equations in 1D setting with similar quadratic nonlinear interactions and Laplacian dissipation. It is therefore considered as a standard benchmark for the analysis of nonlinear advection-diffusion problems. The evolution of the velocity field $u(x, t)$, in a dimensionless form, is given by
\begin{equation}
    \dfrac{\partial u}{\partial t} + u \dfrac{\partial u}{\partial x} = \dfrac{1}{\text{Re}} \dfrac{\partial ^2 u}{\partial x^2}, \label{eq:brg1}
\end{equation}
where $\text{Re}$ is the dimensionless Reynolds number, defined as the ratio of inertial effects to viscous effects. In dimensionless form, the reciprocal of Reynolds number can be denoted as the dimensionless kinematic viscosity $\nu$. Therefore, Eq.~\ref{eq:brg1} can be rewritten as below,
\begin{equation}
    \dfrac{\partial u}{\partial t} + u \dfrac{\partial u}{\partial x} = \nu \dfrac{\partial ^2 u}{\partial x^2}. \label{eq:brg2}
\end{equation}

\subsection{Proper orthogonal decomposition} \label{sec:pod}
The first step for building a projection-based reduced order model is to design a low-order subspace that is capable of capturing the essential features of the system of interest. In the fluid dynamics community, proper orthogonal decomposition (POD) is one of the most popular techniques in this regard. Starting from a collection of system's realizations (called snapshots), POD provides a systematic algorithm to construct a set of orthonormal basis functions (called POD modes) that best describes that collection of snapshot data (in the $\ell_2$ sense). More importantly, those bases are sorted based on their contributions to the system's total energy, making the modal selection a straightforward process. This is a significant advantage compared to other modal decomposition techniques like dynamic mode decomposition, where further sorting and selection criterion has to be carefully defined \cite{bistrian2015improved,bistrian2017method,ahmed2020sampling,kou2017improved}. Usually, the method of snapshots \cite{sirovich1987turbulence} is followed in practice to perform POD efficiently and economically, especially for high dimensional systems. However, we adopt the singular value decomposition (SVD) based approach here for the sake of simplicity and brevity of presentation.

Suppose we have a collection of $N$ system realizations, denoted as $ u(x_i,t_n)$ for $i=1,2,\dots, M$, and $n=1,2,\dots, N$, where $M$ is the number of spatial locations and $N$ is the number of snapshots. Thus, we can build a snapshot matrix $\mathbf{A} \in \mathbb{R}^{M \times N}$ as follows,
\begin{equation}
 \mathbf{A} = \begin{bmatrix}
 u(x_1,t_1) & u(x_1,t_2) & \dots & u(x_1,t_{N}) \\
 u(x_2,t_1) & u(x_2,t_2) & \dots & u(x_2,t_{N}) \\
 \vdots     & \vdots    &  \ddots      & \vdots \\
 u(x_{M},t_1) & u(x_{M},t_2) & \dots & u(x_{M},t_{N}) \\
 \end{bmatrix}.
\end{equation}
Then, a thin (reduced) SVD is performed on $\mathbf{A}$ in the following form, 
\begin{equation} \label{eq:svd1}
    \mathbf{A} = \mathbf{U} \mathbf{\Sigma} \mathbf{V}^T,
\end{equation}
where $\mathbf{U} \in \mathbb{R}^{M \times N}$ is a matrix with orthonormal columns, called the left singular vectors of $\mathbf{A}$ and represent the spatial basis, while the columns of $\mathbf{V} \in \mathbb{R}^{N \times N}$ are the right singular vectors of $\mathbf{A}$, representing the temporal basis. The singular values of $\mathbf{A}$ are stored in descending order as the entries of the diagonal matrix $\mathbf{\Sigma} \in \mathbb{R}^{N \times N}$. Thus, Eq.~\ref{eq:svd1} can be expanded as,
\begin{equation}
 \mathbf{A}  = \begin{bmatrix}
 U_1(x_1) & U_2(x_1) & \dots & U_{N}(x_1) \\
 U_1(x_2) & U_2(x_2) & \dots & U_{N}(x_2) \\
 \vdots     & \vdots    &  \ddots      & \vdots \\
 U_1(x_{M}) & U_2(x_{M}) & \dots & U_{N}(x_{M})
 \end{bmatrix} 
 \begin{bmatrix}
 \sigma_1 &          &        &  \\
          & \sigma_2 &        &  \\
          &          & \ddots &  \\
          &          &        & \sigma_{N}
 \end{bmatrix}
 \begin{bmatrix}
 V_1(t_1) & V_2(t_1) & \dots & V_{N}(t_1) \\
 V_1(t_2) & V_2(t_2) & \dots & V_{N}(t_2) \\
 \vdots     & \vdots    &  \ddots      & \vdots \\
 V_1(t_{N}) & V_2(t_{N}) & \dots & V_{N}(t_{N})
 \end{bmatrix}^{T},
\end{equation}
where $\sigma_1 \ge \sigma_2 \ge \dots \sigma_{N} \ge 0$. For dimensionality reduction purposes, only the first $R$ columns of $\mathbf{U}$, the first $R$ columns of $\mathbf{V}$, and the upper-left $R\times R$ sub-matrix of $\mathbf{\Sigma}$ are considered, corresponding to the largest $R$ singular values. Specifically, the first $R$ columns of $\mathbf{U}$  represent the most effective $R$ POD modes, denoted as $\{\phi_k\}_{k=1}^{R}$ in the rest of the manuscript.  

The velocity field $u(x,t)$ is thus approximated as a linear superposition of the contributions of the first $R$ modes, which can be mathematically expressed as
\begin{equation} \label{eq:uROM}
    u(x,t) = \sum_{k=1}^{R} a_k(t) \phi_k(x),
\end{equation}
where $\phi_k(x)$ are the spatial modes, $a_k(t)$ are the time-dependent modal coefficients (also known as generalized coordinates), and $R$ is the number of retained modes in ROM approximation (i.e., ROM dimension).

\subsection{Galerkin projection} \label{sec:gp}

After constructing a set of POD basis functions, an orthogonal Galerkin projection can be performed to obtain the Galerkin-based ROM (GROM). To do so, the ROM approximation (Eq.~\ref{eq:uROM}) is substituted into the governing equation (Eq.~\ref{eq:brg2}). Noting that the POD bases are only spatial functions (i.e., independent of time) and the modal coefficients are independent of space, we get the following,
\begin{align}
    \left(\sum_{i=1}^{R} \dfrac{\partial a_i}{\partial t}  \phi_i \right) + \left(\sum_{i=1}^{R} a_i \phi_i \right) \left(\sum_{i=1}^{R} a_i \dfrac{\partial \phi_i }{\partial x}  \right) =  \nu  \left(\sum_{i=1}^{R} a_i \dfrac{\partial^2 \phi_i}{\partial x^2}  \right).
\end{align}
We note that the POD basis functions are orthonormal by construction as
\begin{equation}
\langle \phi_i ; \phi_j \rangle = 
     \begin{cases}
       1 &\quad\text{if } i = j \\
       0 &\quad\text{otherwise,}
     \end{cases}    
\end{equation}
where the angle parentheses $\langle \boldsymbol{\bullet} ; \boldsymbol{\bullet} \rangle$ stand for the standard inner product in Euclidean space (i.e., dot product). Then, an inner product with an arbitrary basis function $\phi_k$ can be conducted. Utilizing the orthonormality property of the basis function to simplify ROM derivation, we get the following set of ordinary differential equations (ODEs) representing the tensorial form of GROM,
\begin{align}
    \dfrac{\text{d}a_k}{\text{d}t} &=   \nu \sum_{i=1}^{R} \mathfrak{L}_{i,k} a_i + \sum_{i=1}^{R} \sum_{j=1}^{R} \mathfrak{N}_{i,j,k} a_i a_j, \label{eq:rombrg}
\end{align}
where $\mathfrak{L}$ and $\mathfrak{N}$ are the matrix and tensor of predetermined model coefficients corresponding to linear and nonlinear terms, respectively. Those are precomputed during an offline stage as
\begin{align}
    \mathfrak{L}_{i,k} & = \big\langle \dfrac{\partial^2 \phi_i }{\partial x^2}  ;  \phi_k \big\rangle, \\
    \mathfrak{N}_{i,j,k} &= \big\langle -\phi_i \dfrac{\partial \phi_j}{\partial x};  \phi_k \big\rangle.
\end{align}
Equation~\ref{eq:rombrg} can be rewritten in a compact form as
\begin{equation} \label{eq:romcon}
    \boldsymbol{\dot{\mathbf{a}}} = \mathbf{f}(\mathbf{a}),
\end{equation}
where $\mathbf{a} = [a_1, a_2, \dots, a_R]^T$, and the (continuous-time) model map $\mathbf{f}$ is defined as follows,
\begin{equation*}
\mathbf{f} = 
\begin{bmatrix} 
    \nu \sum_{i=1}^{R} \mathfrak{L}_{i,1} a_i + \sum_{i=1}^{R} \sum_{j=1}^{R} \mathfrak{N}_{i,j,1} a_i a_j \\
    \nu \sum_{i=1}^{R} \mathfrak{L}_{i,2} a_i + \sum_{i=1}^{R} \sum_{j=1}^{R} \mathfrak{N}_{i,j,2} a_i a_j \\
    \vdots \\
    \nu \sum_{i=1}^{R} \mathfrak{L}_{i,R} a_i + \sum_{i=1}^{R} \sum_{j=1}^{R} \mathfrak{N}_{i,j,R} a_i a_j
\end{bmatrix}.
\end{equation*}
Alternatively, Eq.~\ref{eq:romcon} can be used in a discrete-time version as 
\begin{equation} \label{eq:romdis}
    \mathbf{a}^{n+1} = \mathbf{M}(\mathbf{a}^n),
\end{equation}
where $\mathbf{M}(\boldsymbol{\cdot})$ is the discrete-time map obtained by any suitable temporal integration technique. Here, we use the fourth-order Runge-Kutta (RK4) method as follows,
\begin{align}
    \mathbf{a}^{n+1} &= \mathbf{a}^n + \dfrac{\Delta t}{6} (\mathbf{g}_1 + 2\mathbf{g}_2 + 2\mathbf{g}_3 + \mathbf{g}_4),
\end{align}
where 
\begin{align*}
    \mathbf{g}_1 &= \mathbf{f}(\mathbf{a}^n), \\
    \mathbf{g}_2 &= \mathbf{f}(\mathbf{a}^n + \dfrac{\Delta t}{2} \cdot  \mathbf{g}_1), \\
    \mathbf{g}_3 &= \mathbf{f}(\mathbf{a}^n + \dfrac{\Delta t}{2} \cdot  \mathbf{g}_2), \\
    \mathbf{g}_4 &= \mathbf{f}(\mathbf{a}^n + \Delta t \cdot \mathbf{g}_3).
\end{align*}
Thus the discrete-time map defining the transition from time $t_n$ to time $t_{n+1}$ is written as
\begin{align}
\mathbf{M}(\mathbf{a}^n) = \mathbf{a}^n + \dfrac{\Delta t}{6} (\mathbf{g}_1 + 2\mathbf{g}_2 + 2\mathbf{g}_3 + \mathbf{g}_4). 
\end{align}


\section{Long Short-Term Memory Nudging} \label{sec:LSTM}
Due to the quadratic nonlinearity in the governing equation, and consequently the GROM, the \emph{online} computational cost of solving Eq.~\ref{eq:rombrg} is $O(R^3)$ (i.e., it scales cubically with the number of retained modes). Therefore, this has to be kept as low as possible for the feasible implementation of ROM in applications that require near realtime responses (e.g., active control). However, this is often not an easy task for systems with slow decay of the Kolmogorov n-width. Examples include advection-dominated flows with strong nonlinear interactions between wide range of modes. Consequently, the resulting GROM is an intrinsically imperfect model. That is even with the true initial conditions, and absence of numerical errors, the GROM might give inaccurate or false predictions. Indeed,  Carlberg et al. \cite{carlberg2017galerkin} showed that GROM becomes unstable for long time intervals.

Moreover, in most realistic cases, proper specification of the initial state, boundary conditions, and/or model parameters is rarely attainable. This uncertainty in problem definition, in conjunction with model imperfection, poses challenges for accurate predictions. In this study, we put forth a nudging-based methodology that fuses prior model forecast (using imperfect initial condition specification and imperfect model) with the available Eulerian sensor measurements to provide a more accurate posterior prediction. Relating our setting to realistic applications, we build our framework on the assumption that measurements are usually noisy and sparse both in space and time. Nudging has a prestigious history in data assimilation, being a simple and unbiased approach \cite{anthes1974data}. The idea behind nudging is to penalize the dynamical model evolution with the discrepancy between the model's predictions and observations. In other words, the forward model given in Eq.~\ref{eq:romdis} is supplied with a nudging (or correction) term rewritten in the following form, 
\begin{equation}\label{eq:nudge1}
    \mathbf{a}^{n+1} = \mathbf{M}(\mathbf{a}^n) + \mathbf{G}(\mathbf{z}^{n+1}-h(\mathbf{a}^{n+1})),
\end{equation}
where $\mathbf{G}$ is called the nudging (gain) matrix and $\mathbf{z}$ is the set of measurements (observations), while $h(\cdot)$ is a mapping from model space to observation space. For example, $h(\cdot)$ can be a reconstruction map, from ROM space to FOM space. In other words, $h(\mathbf{a})$ represents the ``model prediction of the measured quantity'', while $\mathbf{z}$ is the ``actual'' observations. Given the simplicity of Eq.~\ref{eq:nudge1}, the specification/definition of the gain matrix $\mathbf{G}$ is not as straightforward \cite{zou1992optimal,vidard2003determination,auroux2005back,lakshmivarahan2013nudging}.

In the proposed framework, we utilize recurrent neural networks, namely the long short-term memory (LSTM) variant, to define a nudging map. In particular, Eq.~\ref{eq:nudge1} implies that each component of $\mathbf{a}^{n+1}$ (i.e., $a_1,a_2\dots,a_R$) is corrected using a \emph{linear} superposition of the the components of $\mathbf{z}^{n+1}-h(\mathbf{a}^{n+1})$, weighted by the gain matrix. Here, we relax this linearity assumption and generalize it to a possibly nonlinear mapping $\mathbf{C}(\mathbf{a}, \mathbf{z})$ as,
\begin{equation}\label{eq:nudge2}
    \mathbf{a}^{n+1} = \mathbf{M}(\mathbf{a}^n) + \mathbf{C}(\mathbf{a}_b^{n+1}, \mathbf{z}^{n+1}),
\end{equation}
where the map $\mathbf{C}(\mathbf{a}, \mathbf{z})$ is learnt (or fit) using an LSTM neural network, and $\mathbf{a}_b^{n+1}$ is the prior model prediction computed using imperfect model and/or imperfect initial conditions (called background in data assimilation terminology), defined as $\mathbf{a}_b^{n+1} = \mathbf{M}(\mathbf{a}^n)$. Thus, Eq.~\ref{eq:nudge2} can be rewritten as follows,
\begin{equation}\label{eq:nudge3}
    \mathbf{a}^{n+1} = \mathbf{a}_b^{n+1} + \mathbf{C}(\mathbf{a}_b^{n+1}, \mathbf{z}^{n+1}).
\end{equation}

In order to learn the map $\mathbf{C}(\mathbf{a}_b,\mathbf{z})$, we consider the case with an imperfect model, defective initial conditions, and noisy observations. Moreover, we suppose sensors are sparse in space, and measurement signals are sparse in time too. Specifically, we use sensors located at a few equally-spaced grid points, but a generalization to off-grid sensor placement is possible. Also, we assume sensors send measurement signals every $\tau$ time units. To mimic sensor measurements and noisy initial conditions, we run a twin experiment as follows,
\begin{enumerate}
    \item Solve the FOM equation (i.e., Eq.~\ref{eq:brg2}) and sample \emph{true} field data ($u_{true}(x,t_n)$) each $\tau$ time units. In other words, store $u_{true}(x,t_n)$ at $t_n\in \{0,\tau,2\tau,\dots T\}$ where $T$ is the total (maximum) time and $\tau$ is the time window over which measurements are collected.
    \item Define erroneous initial field estimate as $u_{err} (x,t_n) = u_{true}(x,t_n) + \epsilon_b$, where $t_n\in \{0,\tau,2\tau,\dots T-\tau\}$. $\epsilon_b$ stands for noise in initial state estimate, assumed as white Gaussian noise with zero mean and covariance matrix $B$ (i.e., $\epsilon_b \sim {\cal{N}}(0,B)$).
    \item Define sparse and noisy measurements as $\mathbf{z} = u_{true}(x_{Obs},t_n) + \epsilon_m $, for $t_n\in \{\tau,2\tau,\dots T\}$. Similarly, $\epsilon_m$ stands for the measurements noise, assumed to be white Gaussian noise with zero mean and covariance matrix $Q$ (i.e., $\epsilon_m \sim {\cal{N}}(0,Q)$).
\end{enumerate}
For LSTM training data, we project the erroneous field estimates (from Step 2) onto the POD basis functions to get the erroneous POD modal coefficients (i.e., $\mathbf{a}_{err}(t_n)$, for $t_n\in \{0,\tau,2\tau,\dots T-\tau\}$. Then, we integrate those erroneous coefficients for $\tau$ time units to get the background prediction $\mathbf{a}_b(t_n)$, for $t_n\in \{\tau,2\tau,\dots T\}$. 

Then, we train the LSTM using $\mathbf{a}_b(t_n)$ and $\mathbf{z}(t_n)$ as inputs, and set the target as the correction $(\mathbf{a}_{true}(t_n) - \mathbf{a}_b(t_n))$, for $t_n\in \{\tau,2\tau,\dots T\}$. The true modal coefficients ($\mathbf{a}_{true}$) are obtained by projecting the \emph{true} field data (from Step 1) onto the POD bases, where the projection is defined via the inner product as $a_k(t) = \langle u(x,t) ; \phi_k(x) \rangle$. In order to enrich the training data set, Step 2 and Step 3 are repeated several times giving an ensemble of erroneous state estimates and noisy measurements at every time instant of interest. Each member of those ensembles represents one training sample. This also assists the LSTM network to handle wider range of noise.

We emphasize that the proposed LSTM-Nudge approach not only cures model imperfection (i.e., provides model closure and accounts for any missing physical processes) but also treats uncertainties in initial state estimates. Moreover, the field measurements (i.e., the nudging data) are assumed to be sparse and noisy to mimic real-life situations. 

\section{Results} \label{sec:res}
\begin{figure}[ht]
\centering
\includegraphics[trim= 0 0 0 0, clip, width=0.55\textwidth]{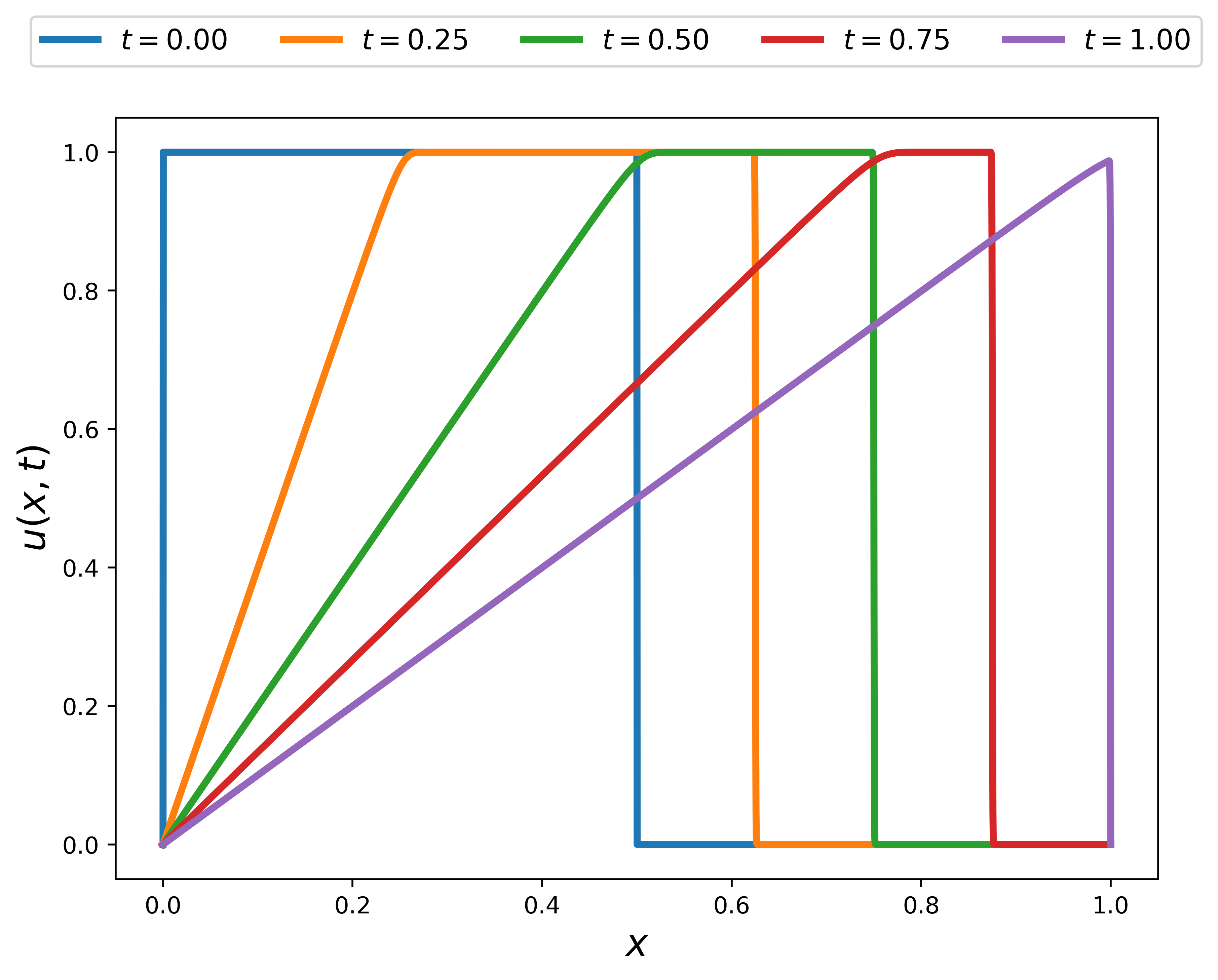}
\caption{Evolution of the FOM velocity field for the 1D Burgers problem, characterized by a moving shock with initial square wave.}
\label{fig:FOM}
\end{figure}

We test the proposed methodology using the 1D Burgers problem introduced in Sec.~\ref{sec:gov}. In particular, we consider a domain of dimensionless length of one, with a square wave as initial condition defined as,

\begin{equation}\label{eq:IC}
u(x,0) = \begin{cases}
            1 , &\quad\text{if } 0 < x \le 0.5 \\
            0 , &\quad\text{if } 0.5 < x \le 1.0,
        \end{cases}
\end{equation}
with zero Dirichlet boundary conditions, $u(0,t) = u(1,t) = 0$. We solve Eq.~\ref{eq:brg1} at $\text{Re} = 10^4$ for $t \in [0,1]$. For numerical computations, we use a family of fourth order compact schemes for spatial derivatives \cite{lele1992compact}, and skew-symmetric formulation for the nonlinear term. For the FOM simulation, we use a time step of $10^{-4}$ over a spatial grid of $4096$, and for POD basis generation, we collect 100 snapshots (i.e., every 100 FOM time steps). The temporal evolution of the 1D Burgers problem using the described setup is shown in Fig.~\ref{fig:FOM}, illustrating the advection of the shock wave over time.

For ROM computations, 6 modes are retained in the reduced order approximation (i.e.,  $R=6$) and a time step of $0.01$ is adopted for the temporal integration of GROM equations. In order to implement the LSTM-Nudge approach, we begin at erroneous initial condition defined as $u_{err}(x,0) = u_{true}(x,0) + \epsilon_b$, where $u_{true}(x,0)$ is defined with Eq.~\ref{eq:IC}, and $\epsilon_b$ is a white Gaussian noise with zero mean and covariance matrix $B$. For simplicity, we assume $B=\sigma_b^2 \mathbf{I}$, where $\sigma_b$ is the standard deviation in ``background'' estimate of the initial condition and $\mathbf{I}$ is the identity matrix. We note that this formulation implies that our estimates of the initial velocity field at different spatial locations are uncorrelated. As nudging field data, we locate sensors to measure the velocity field $u(x,t)$ every 256 grid points (i.e., a total of 17 sensors with $s_{freq} = 256$, where $s_{freq}$ is the number of spatial steps between sensors locations), and collect measurements every 10 time steps (i.e., each $0.1$ time unit with $t_{freq}=10$, where $t_{freq}$ is the number of time steps between measurement signals). To account for noisy observations, white Gaussian noise of zero mean and covariance matrix of $Q$ is added to the true velocity field obtained from the FOM simulation. Similar to $B$, we set $Q=\sigma_m^2 \mathbf{I}$, where $\sigma_m$ is the standard deviation of measurement noise. This assumes that sensors measurements are not correlated to each other, and all sensors have similar quality (i.e., add similar amounts of noise to the measurements). As a base case, we set $\sigma_b=1$, and $\sigma_m=1$.

\begin{figure}[H]
\centering
\includegraphics[trim= 0 0 0 0, clip, width=0.95\textwidth]{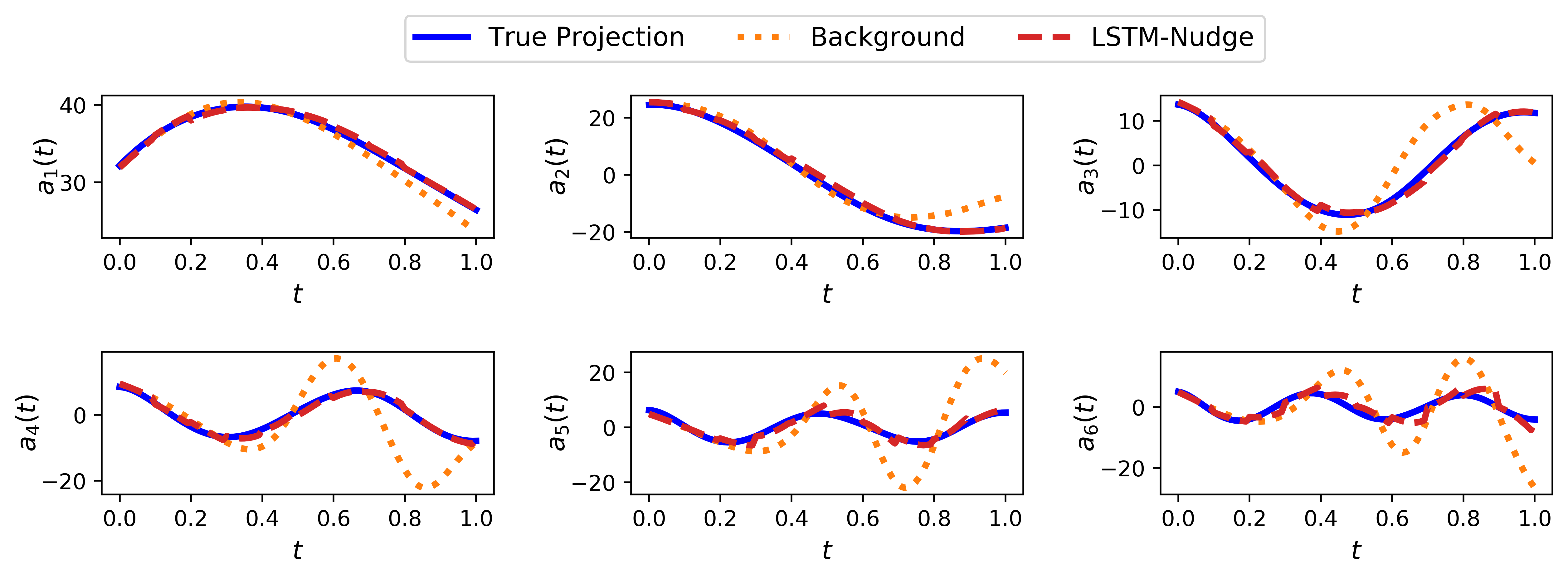}
\caption{Temporal evolution of the POD modal coefficients for the 1D Burgers problem.}
\label{fig:a}
\end{figure}

The procedure presented in Sec.~\ref{sec:LSTM} is applied using the numerical setup described above, and compared against the reference case of GROM with the erroneous initial condition and inherent model imperfections due to modal truncation. In Fig.~\ref{fig:a}, the temporal evolution of the POD modal coefficients is shown for the true projection, background, and LSTM-Nudge results. The true projection results are obtained by the projection (i.e., via inner product) of the true FOM field at different time instants onto the corresponding basis functions. The background trajectory is the reference solution obtained by standard GROM using the erroneous initial condition, without any closure or corrections. It can be seen that the background trajectory gets off the true trajectory by time as a manifestation of model imperfection. Also, note that the background solution does not begin from the same point as true projection due to the noise in initial condition. On the other hand, the LSTM-Nudge predictions almost perfectly match the true projection solution, implying that the approach is capable of blending noisy observations with a prior estimate to gain more accurate predictions.

\begin{figure}[H]
\centering
\includegraphics[trim= 0 0 0 0, clip, width=0.95\textwidth]{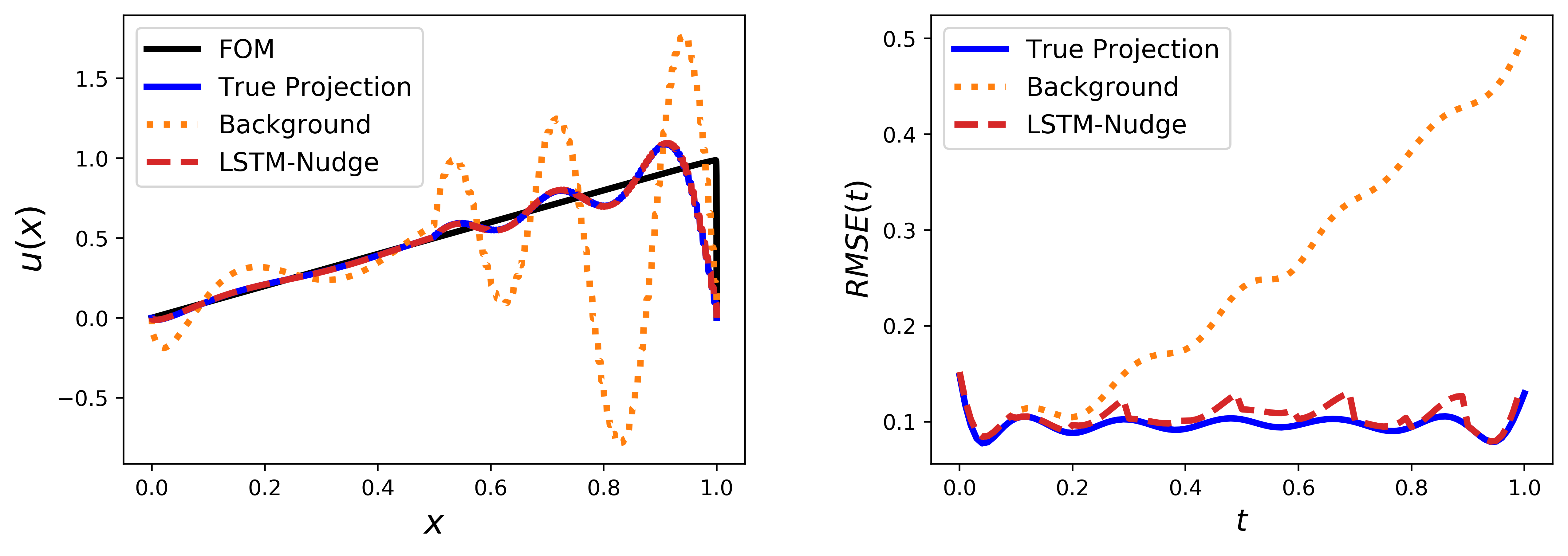}
\caption{Final velocity field (at $t=1$) [left] and the root mean-squares error [right] for the 1D Burgers problem.}
\label{fig:u}
\end{figure}

In order to better visualize the predictive capabilities of the LSTM-Nudge methodology, we compute the reconstructed velocity field using Eq.~\ref{eq:uROM}. Moreover, the root mean-squares error ($RMSE$) of the reconstructed field with respect to the FOM solution is calculated as a function of time as follows,

\begin{equation} \label{eq:RMSE}
    RMSE(t) = \sqrt{\dfrac{1}{M}\sum_{i=1}^{M}{\bigg(u_{FOM}(x_i,t) - u_{ROM}(x_i,t) \bigg)^2} },
\end{equation}
where $u_{FOM}$ is the true velocity field obtained from solving the FOM equation, while $u_{ROM}$ is the reduced order approximation computed through true projection, background (reference) solution, or LSTM-Nudge method. The reconstructed velocity field at final time (i.e., at $t=1$) is shown in Fig.~\ref{fig:u} along with the $RMSE$ as a function of time. As described before, the true projection solution is simply the projection of the FOM field onto the reduced POD space, and it represents the optimal solution that can be approximated using a linear subspace spanned by $R$ modes. In order to get rid of those Gibbs oscillations, we would need either a larger number of modes or a more representative subspace (e.g., through partitioning or auto-encoders).  Therefore, it is fair to compare our results against the true projection solution, rather than the FOM since we do not address any issues regarding the resolution or representability capabilities of the POD subspace.

\subsection{Effect of noise} \label{sec:noise}

\begin{figure}[H]	
	\centering
	\begin{subfigure}[t]{0.48\linewidth}
		\centering
		\includegraphics[trim= 390 0 0 0, clip, width=1\linewidth]{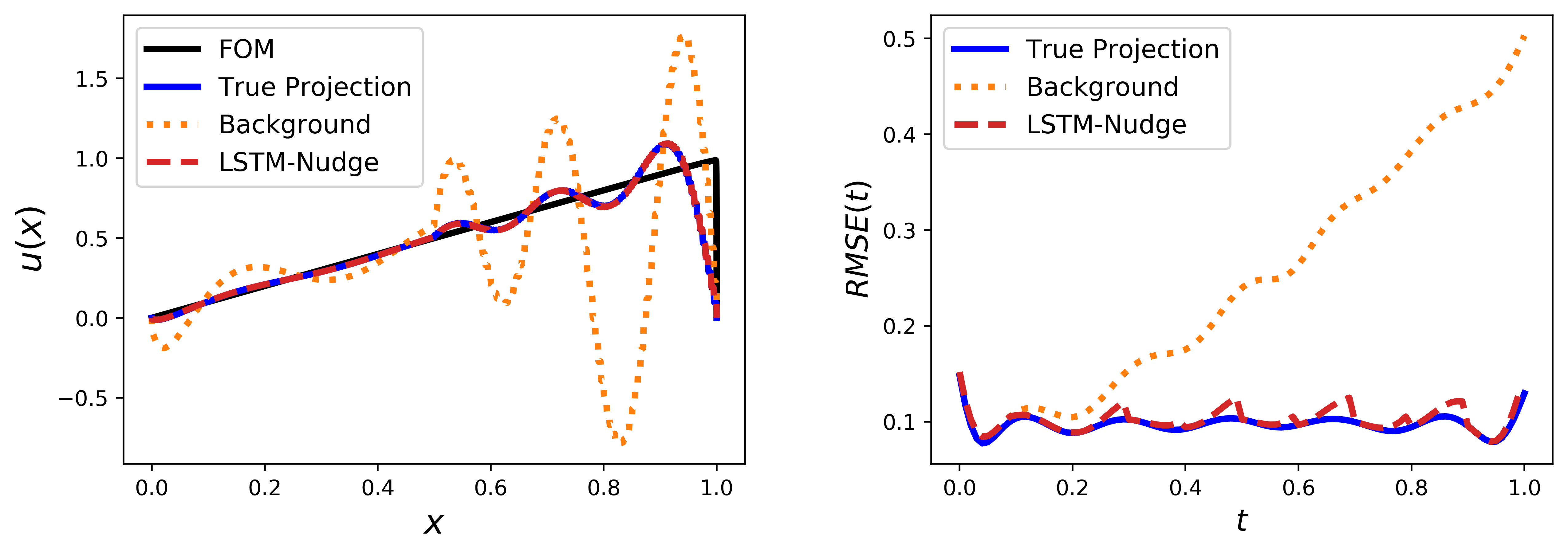}
		\caption{$\sigma_m = 0.1, \sigma_b = 1.0$}\label{fig:noisea}		
	\end{subfigure}
	\quad
	\begin{subfigure}[t]{0.48\linewidth}
		\centering
		\includegraphics[trim= 390 0 0 0, clip, width=1\linewidth]{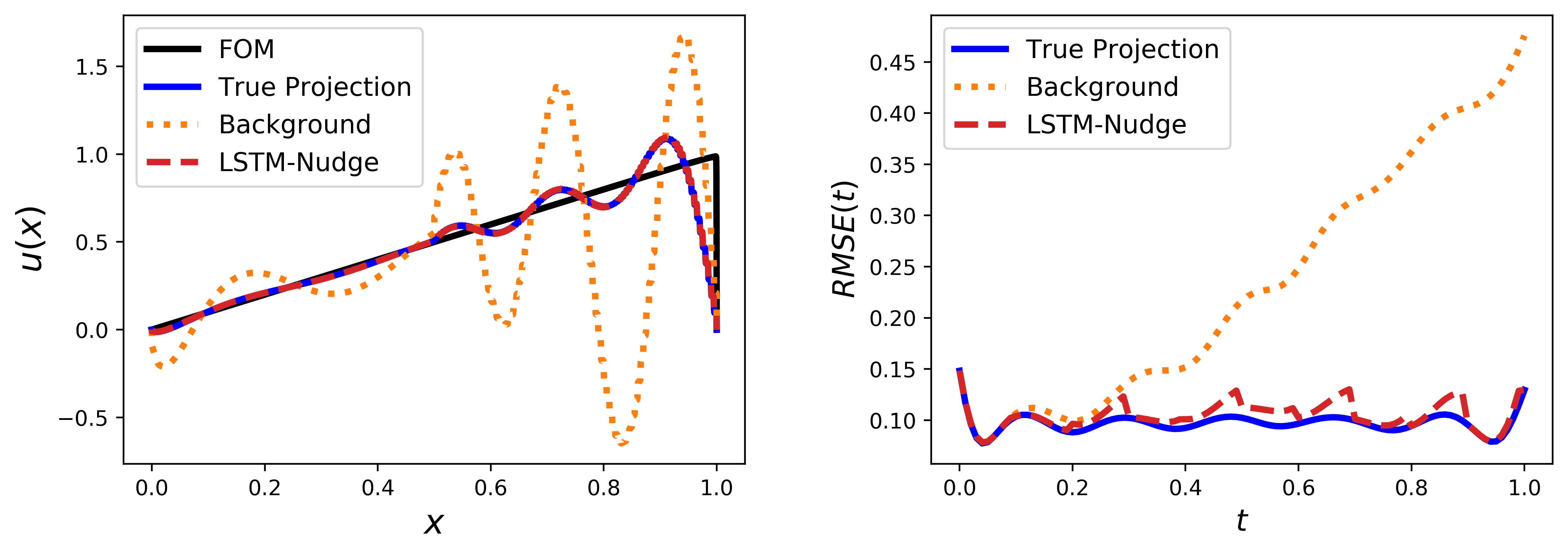}
		\caption{$\sigma_m = 1.0, \sigma_b = 0.1$}\label{fig:noiseb}		
	\end{subfigure}
    \\ \bigskip
    \begin{subfigure}[t]{0.48\linewidth}
		\centering
		\includegraphics[trim= 390 0 0 0, clip, width=1\linewidth]{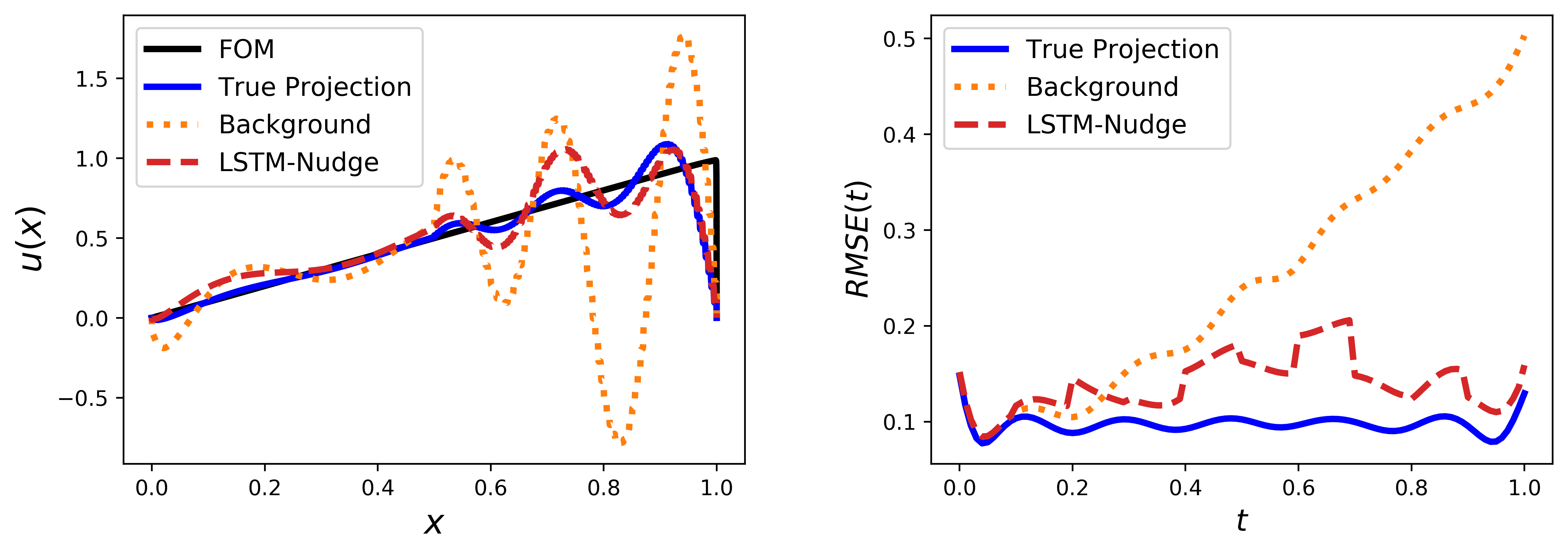}
		\caption{$\sigma_m = 10.0, \sigma_b = 1.0$}\label{fig:noisec}		
	\end{subfigure}
	\quad
	\begin{subfigure}[t]{0.48\linewidth}
		\centering
		\includegraphics[trim= 390 0 0 0, clip, width=1\linewidth]{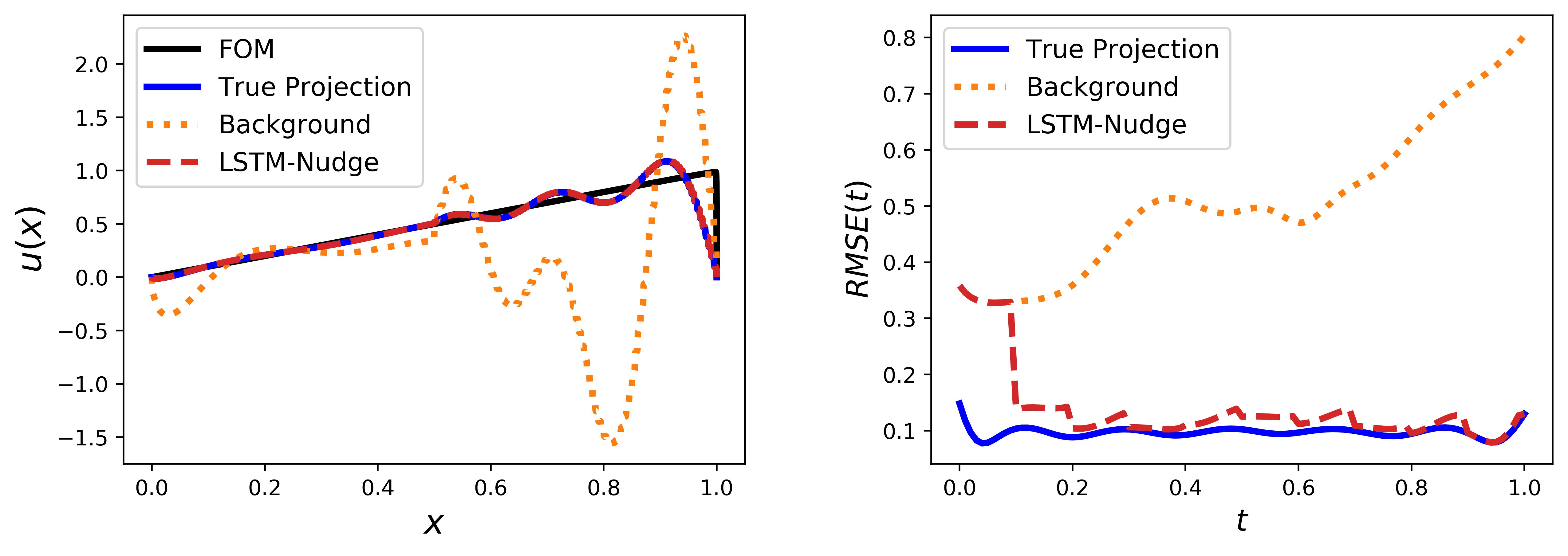}
		\caption{$\sigma_m = 1.0, \sigma_b = 10.0$}\label{fig:noised}		
	\end{subfigure}
	  \\ \bigskip
    \begin{subfigure}[t]{0.48\linewidth}
		\centering
		\includegraphics[trim= 390 0 0 0, clip, width=1\linewidth]{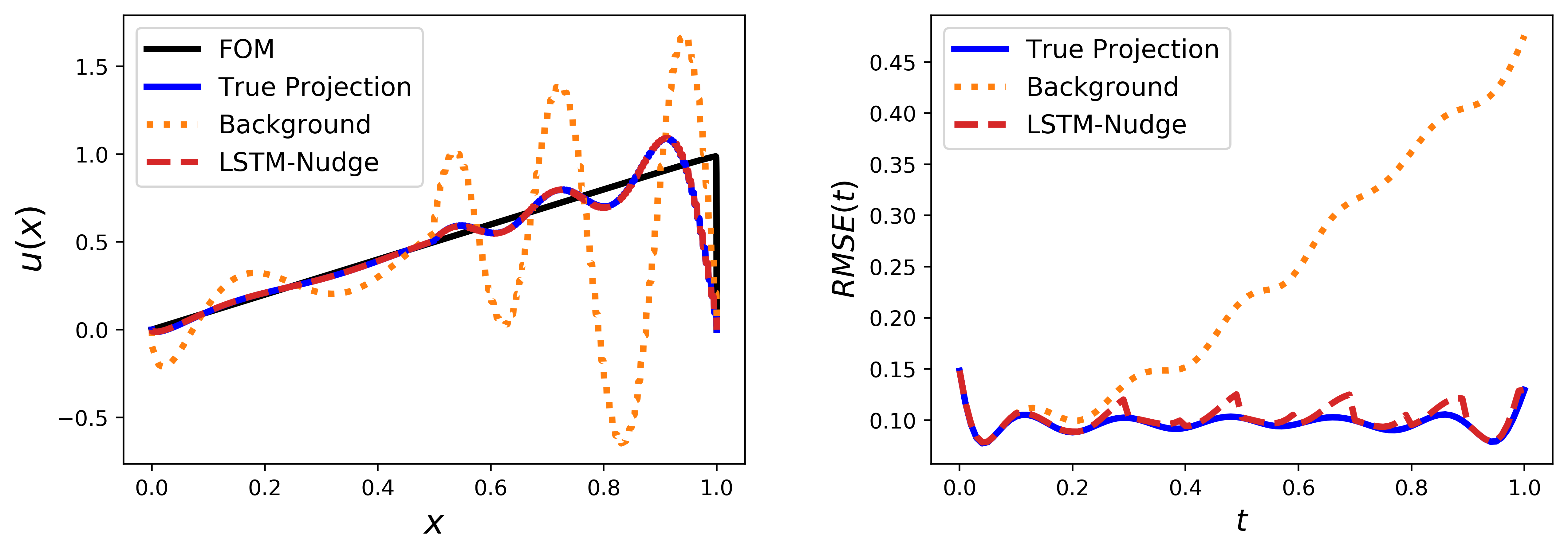}
		\caption{$\sigma_m = 0.1, \sigma_b = 0.1$}\label{fig:noisee}		
	\end{subfigure}
	\quad
	\begin{subfigure}[t]{0.48\linewidth}
		\centering
		\includegraphics[trim= 390 0 0 0, clip, width=1\linewidth]{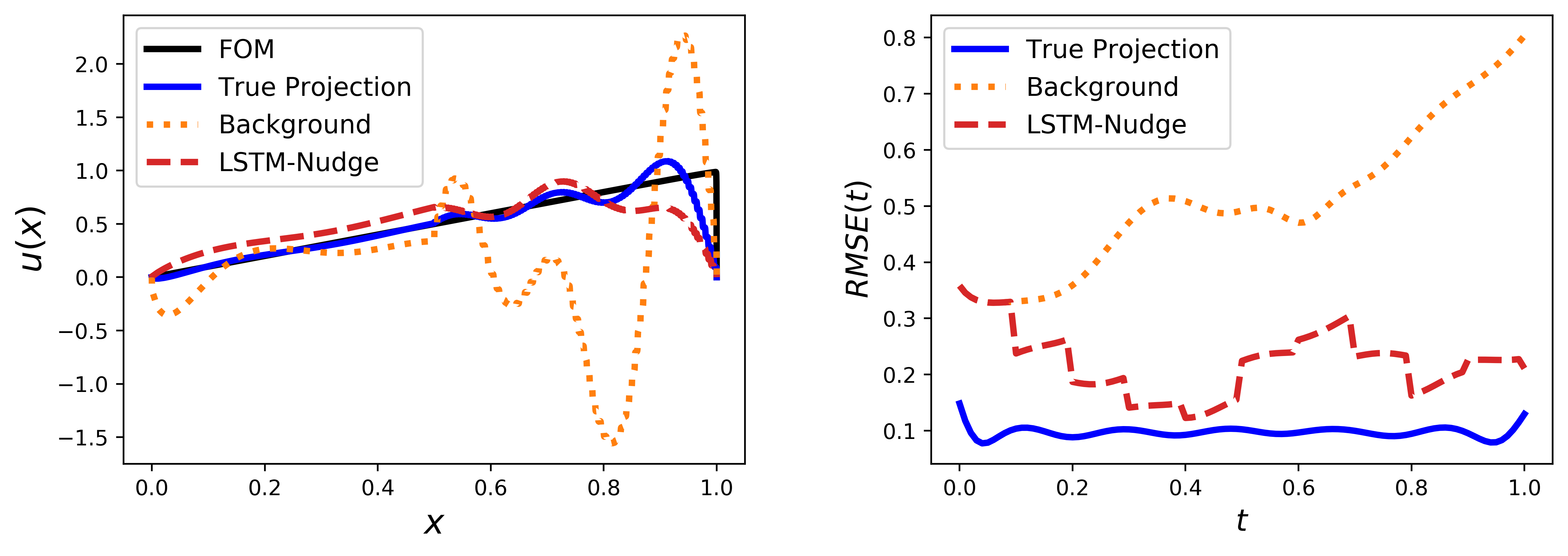}
		\caption{$\sigma_m = 10.0, \sigma_b = 10.0$}\label{fig:noisef}		
	\end{subfigure}
	\caption{Root mean-squares error in reconstructed velocity field, with different levels of background and measurement noises.}\label{fig:noise}
\end{figure}

\begin{figure}[H]
\centering
\includegraphics[trim= 75 0 45 0, clip, width=\textwidth]{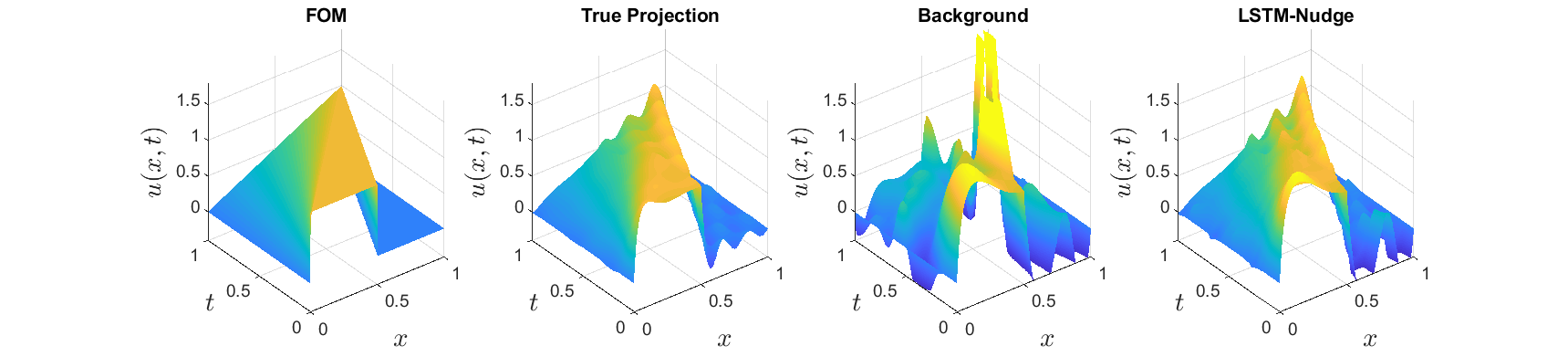}
\caption{Surface plots for the spatio-temporal evolution of the 1D Burgers problem with $\sigma_b=10$, and $\sigma_m=1$.}
\label{fig:surf1}
\end{figure}

Here, we investigate the effect of noise (both in initial condition and measurements) on the performance of the LSTM-Nudge framework. In other words, we study how much noise it can handle sufficiently. For the training phase, the LSTM was trained using noisy data with $\sigma_b = \sigma_m = 1.0$. Now, we test using data with smaller and larger amounts of noise. In particular, we vary $\sigma_b$ and $\sigma_m$ between $0.1$, $1.0$, and $10.0$. Readers should be aware that the true velocity field spans between $0$ and $1$. Thus, a noise with a standard deviation of 10 is an extreme case, corresponding to very cheap sensors. In Fig.~\ref{fig:noise}, we show the root mean-squares error of the reconstructed velocity fields based on true projection, background solution, and LSTM-Nudge predictions using different levels of noise. We can see that the LSTM-Nudge is performing very well, compared to the background solution, and almost matching the true projection results. More importantly, we find that the prediction accuracy is more dependent on measurement noise than noise in the initial condition. For instance, the LSTM-Nudge almost recovers the true state estimate very shortly using adequate measurements (i.e., Fig.~\ref{fig:noised}). This is even more visualized in the surface plots in Fig.~\ref{fig:surf1}. In contrast, the model imperfections cannot be cured well at many time instances using highly noisy observations, even with moderate noise in initial condition (e.g., see Fig.~\ref{fig:noisec}). The situation is worse in Fig.~\ref{fig:noisef}, with severe noise in both initial conditions and measurements. 

The dependence of the framework on measurements noise significantly more than background noise might be attributed to the input features we are using in the LSTM architecture. It includes a combination between the background modal coefficients (obtained from erroneous initial condition) and the observed velocity field at sparse locations. For the modal coefficients, the erroneous initial conditions are first projected onto the POD subspace to obtain the initial (erroneous) modal coefficients, which are then integrated in time for $\tau$ time units. This projection is known to filter a large amount of added noise, which can be seen as a preprocessing step to reduce the effect of initial condition perturbation. That is why we can barely see a difference between erroneous trajectory and true projection at time zero, except for the extreme case with $\sigma_b=10$. On the other hand, observations are fed to the LSTM network as is, without any preprocessing or prior treatment. This makes the predictions more sensitive to the measurement quality.

\subsection{Effect of measurements sparsity} \label{sec:sparse}
\begin{figure}[H]
\centering
\includegraphics[trim= 0 0 0 0, clip, width=0.95\textwidth]{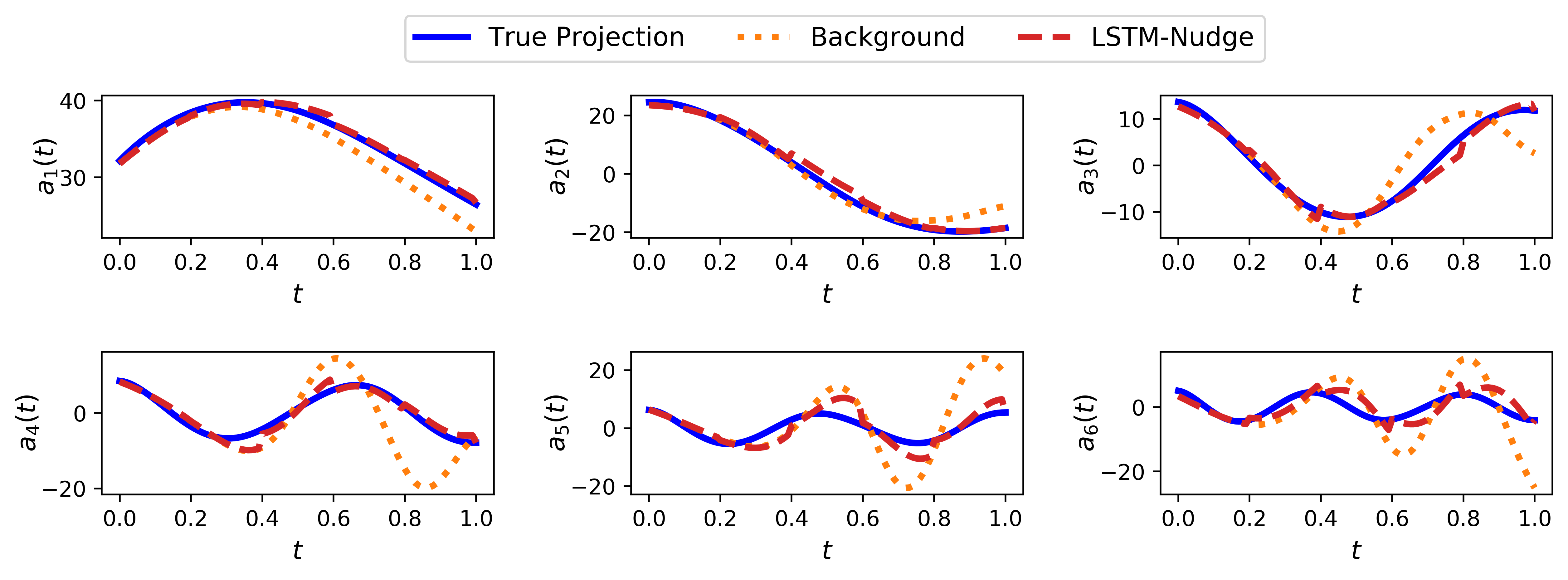}
\caption{Temporal evolution of the POD modal coefficients for the 1D Burgers problem, when measurements are taken every 20 time steps.}
\label{fig:tsparsea}
\end{figure}

Since the LSTM-Nudge in ROM context is found to be relatively sensitive to sensor quality (observational noise), we study the effect of measurement sparsity as well. In particular, both the temporal sparsity (i.e., frequency of measurement signals) and spatial sparsity (i.e., number of sensors) are explored. For the base case, we are collecting measurements every 10 time steps (i.e., $t_{freq}=10, \tau=0.1$), and sensors are placed every 256 grid points (i.e., $s_{freq}=256$). First, we consider the case when measurement signals are only available every 20 time steps (i.e., $t_{freq}=20$), with the same spatial sparsity (i.e., 17 sensors). The POD modal coefficients are plotted in Fig.~\ref{fig:tsparsea} as predicted by LSTM-Nudge, compared to the background case (without corrections) and true projection trajectory. We find that the framework is sufficiently able to handle this variation in measurement signal frequency. We note here that we use the same LSTM network, trained using the base case data (i.e., trained with $t_{freq}=10$ and tested for $t_{freq}=20$). For all, we use the same level of noise as before (i.e., $\sigma_b=\sigma_m=1$). 

\begin{figure}[H]
\centering
\includegraphics[trim= 0 0 0 0, clip, width=0.95\textwidth]{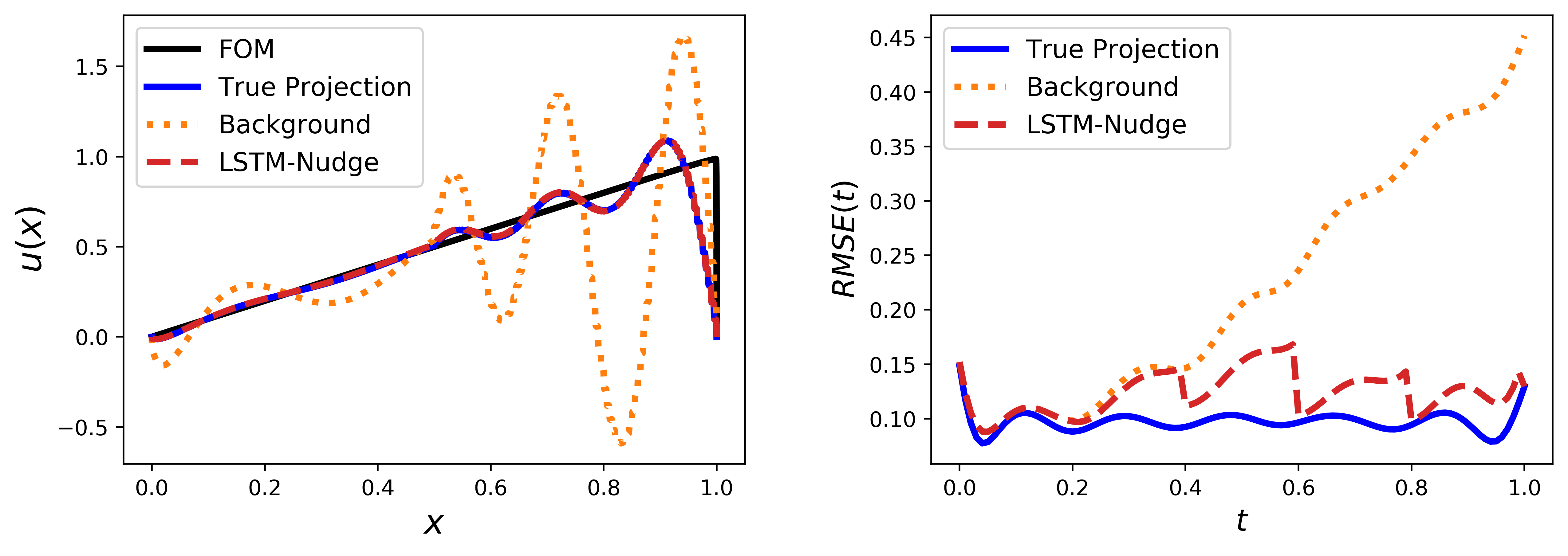}
\caption{Final velocity field (at $t=1$) [left] as well as the root mean-squares error [right] for the 1D Burgers problem, when measurements are taken every 20 time steps.}
\label{fig:tsparseu}
\end{figure}

We also plot the final velocity field reconstructed through the background solution, and LSTM-Nudge compared to both FOM and true projection in Fig.~\ref{fig:tsparseu} as well as the root mean-squares errors at different times. Although we can see a discrepancy between true projection and LSTM-Nudge at several time instants, we notice a jump in LSTM-Nudge solution towards the true projection almost every $0.2$ time units. This is consistent with the fact that LSTM-Nudge goes into effect every 20 time steps, when measurements become available, implying that the LSTM-Nudge is still capable of sufficiently rectifying the model and state imperfections whenever measurements are received.

\begin{figure}[ht]	
	\centering
	\begin{subfigure}[t]{0.48\linewidth}
		\centering
		\includegraphics[trim= 390 0 0 0, clip, width=1\linewidth]{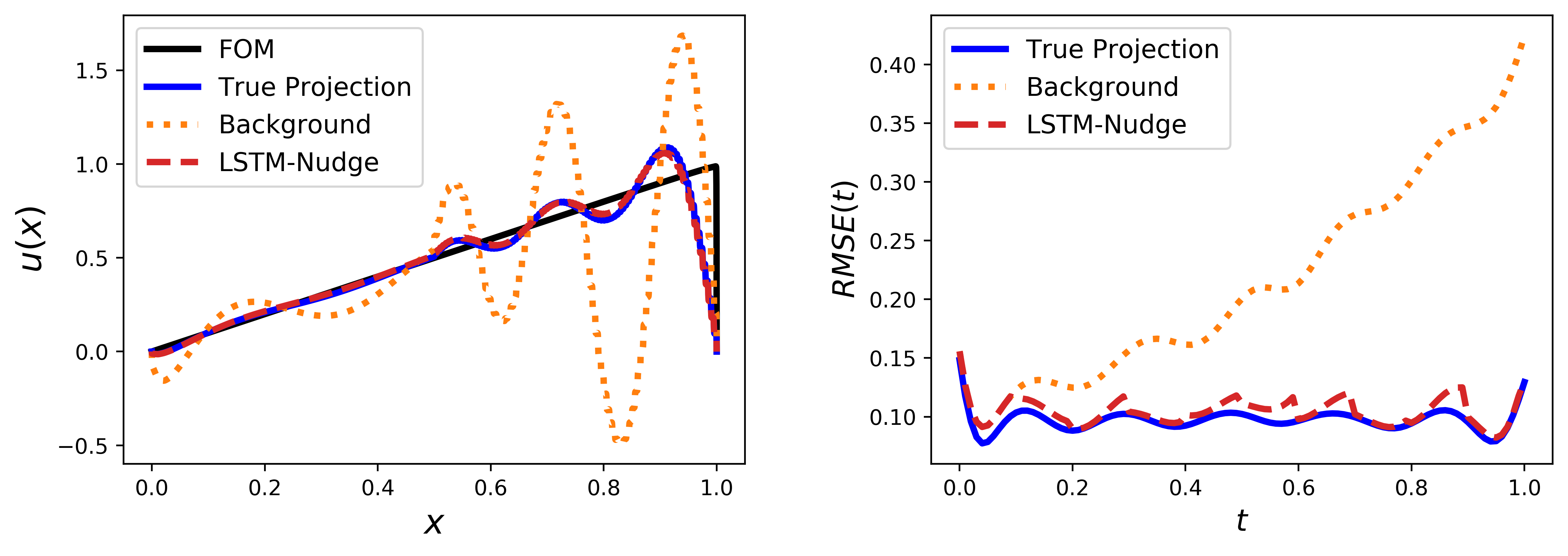}
		\caption{sensors located every 128 grid point}\label{fig:ssp1}		
	\end{subfigure}
	\quad
    \begin{subfigure}[t]{0.48\linewidth}
		\centering
		\includegraphics[trim= 390 0 0 0, clip, width=1\linewidth]{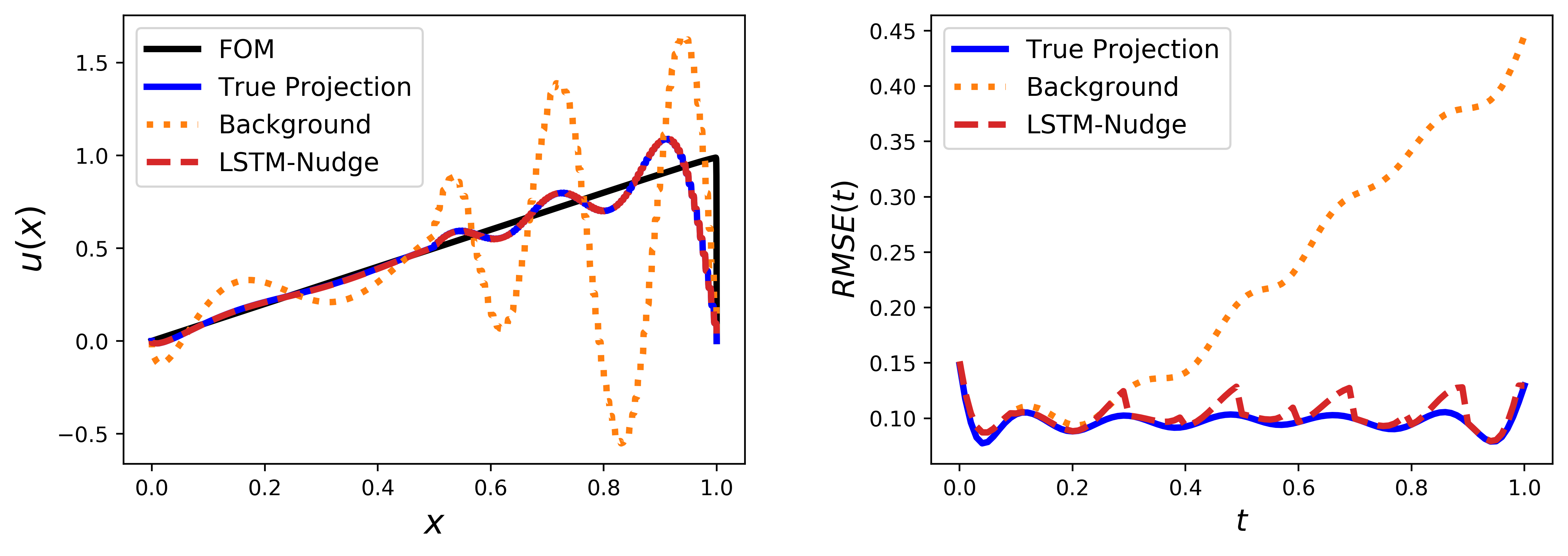}
		\caption{sensors located every 512 grid point}\label{fig:ssp2}		
	\end{subfigure}
	\\ \bigskip
	\begin{subfigure}[t]{0.48\linewidth}
		\centering
		\includegraphics[trim= 390 0 0 0, clip, width=1\linewidth]{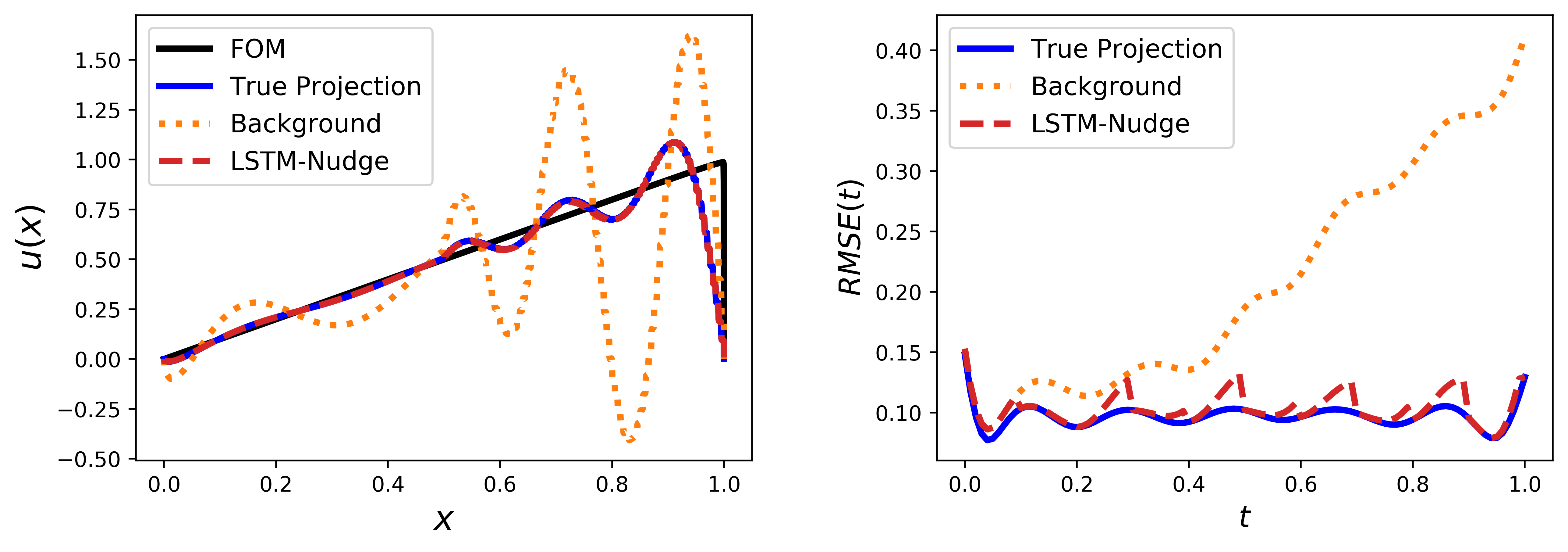}
		\caption{sensors located every 1024 grid point}\label{fig:ssp3}		
	\end{subfigure}
	\quad
    \begin{subfigure}[t]{0.48\linewidth}
		\centering
		\includegraphics[trim= 390 0 0 0, clip, width=1\linewidth]{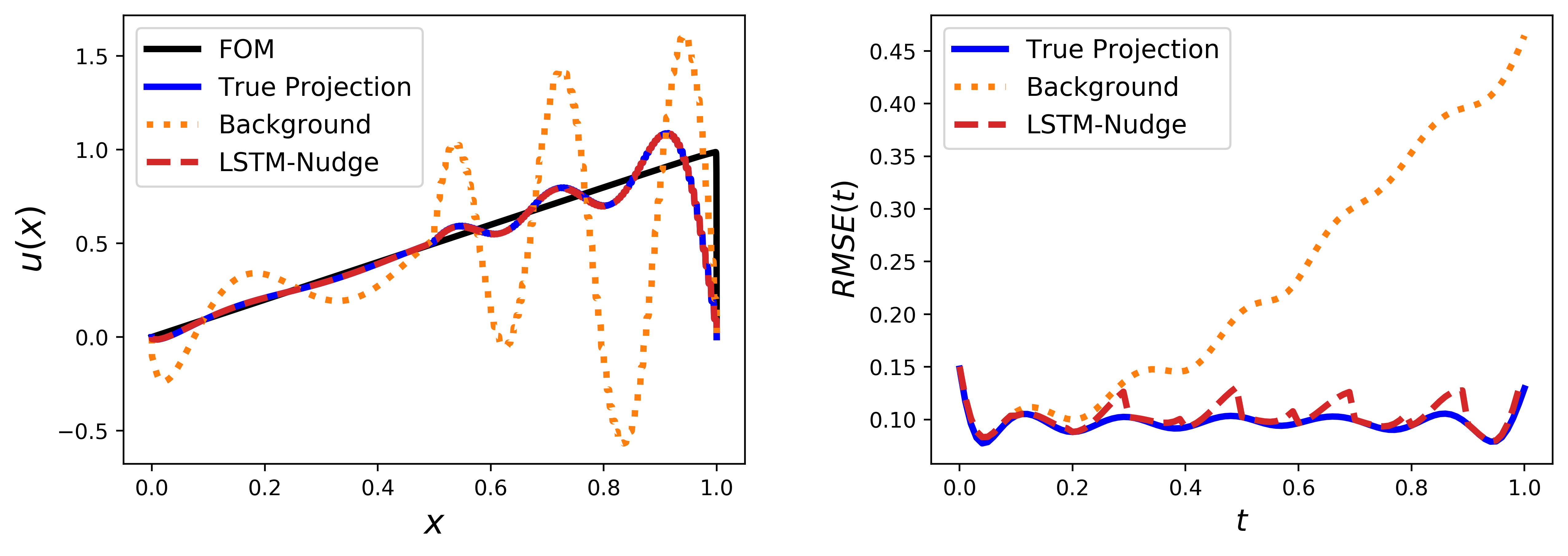}
		\caption{sensors located every 2048 grid point}\label{fig:ssp4}		
	\end{subfigure}
	
	\caption{Root mean-squares error in reconstructed velocity field, with different number of sensors located sparsely at grid points. } \label{fig:ssparse}
\end{figure}

In order to examine the effect of spatial sparsity (number of sensors), we vary the spatial frequency (i.e., number of grid points between sensors) as $s_{freq} \in \{128, 512, 1024, 2048\}$. The first case (i.e., $s_{freq}=128$) corresponds to more sensors than base case, while the others correspond to more sparse measurement points (less sensors). From Fig.~\ref{fig:ssparse}, we can deduce that the effect of number of sensors is minimal in this case, even with very few sensors (e.g., 3 sensors in Fig~\ref{fig:ssp4}). However, we should state here that each of those cases requires retraining the LSTM network with the relevant number of measurement points. This is because the LSTM for the base case has an input dimension of 23 (i.e., 6 modal coefficients and 17 measurements), and changing the number of measurements would require a different size of input vector. Although we assume equally-spaced and collocated sensors (i.e., placed exactly on the numerical grid), compressive sensing ideas can be adopted to intelligently locate sensors for optimal performance.

\subsection{Effect of measured quantity} \label{sec:measure}
\begin{figure}[H]
\centering
\includegraphics[trim= 0 0 0 0, clip, width=0.95\textwidth]{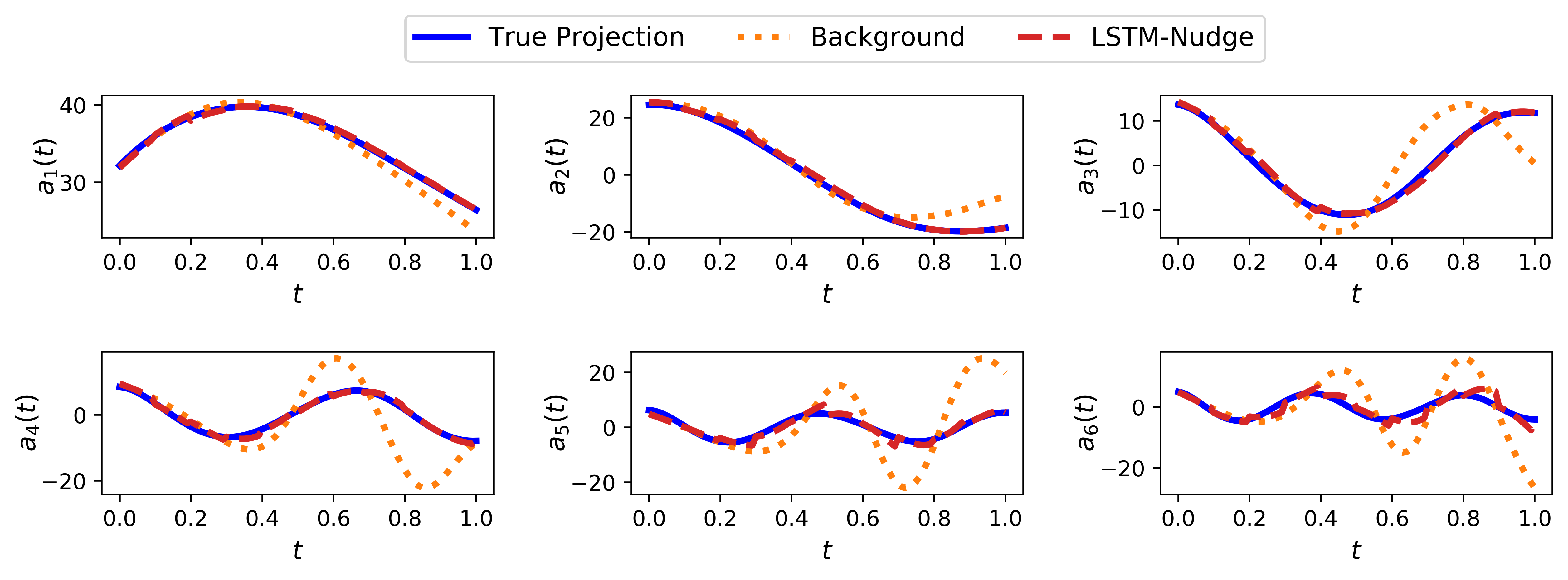}
\caption{Temporal evolution of the POD modal coefficients for the 1D Burgers problem, with $u^2$ as the available measurements.}
\label{fig:au2}
\end{figure}

\begin{figure}[H]
\centering
\includegraphics[trim= 0 0 0 0, clip, width=0.95\textwidth]{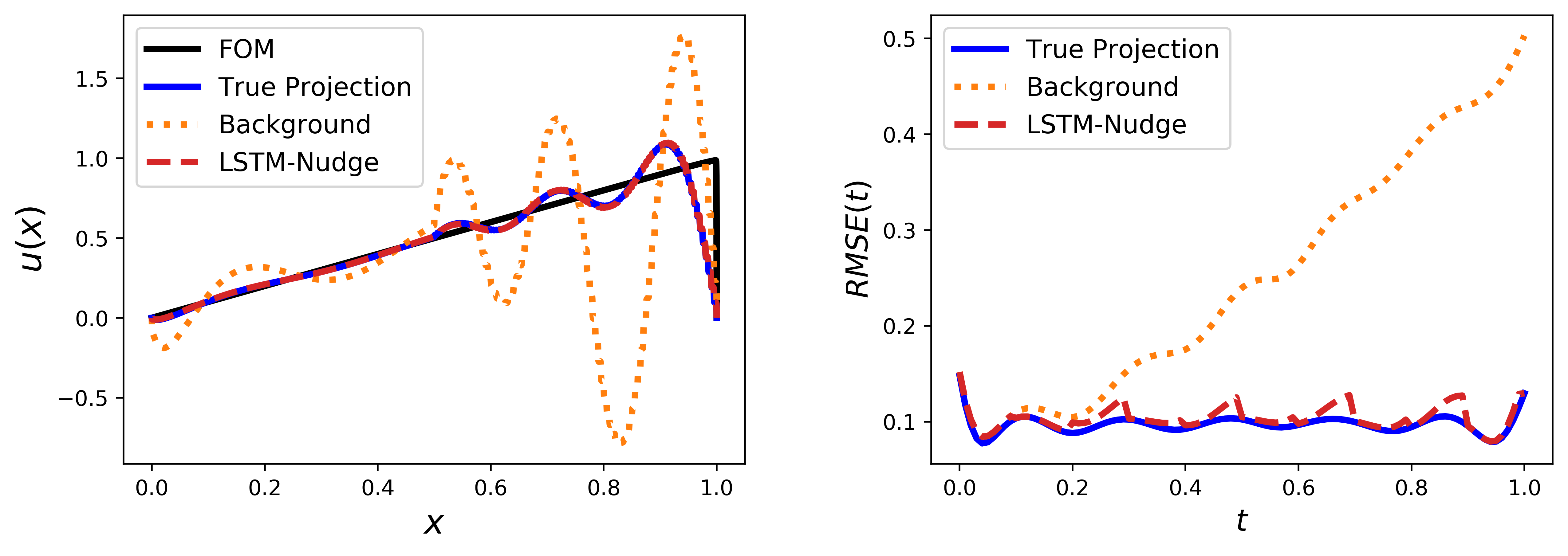}
\caption{Final velocity field (at $t=1$) [left] as well as the root mean-squares error [right] for the 1D Burgers problem,  with $u^2$ as the available measurements.}
\label{fig:uu2}
\end{figure}

As described in Sec.~\ref{sec:LSTM}, the input to the LSTM-Nudge framework is a combination of modal coefficients (i.e., GROM state variables), and direct measurements (i.e., velocity field) without constraining any mapping between them. The LSTM has shown brilliant effectiveness learning the map between model state and observations to approximate the required correction/nudging. In this section, we elaborate more on this feature by exploring the performance of the LSTM-Nudge with a different measured quantity. In practice, direct field variable measurement may not be feasible. For example, the dynamics of sea surface
temperature can be only inferred from satellite measurement of radiated thermal energy. Defining a map between observable quantity and model (state) variable is not usually straightforward. Hence, utilizing neural networks strengths of discovering underlying patterns and relations to learn such maps is highly desired. Instead of measuring the velocity $u(x,t)$, we hypothesize that we can only observe the square of this velocity field (i.e., $u^2(x,t)$). This is related to the kinetic energy of the flow. We repeat the LSTM training using the new input features (i.e., modal coefficients and square of velocity), and test using the base case parameters (i.e., $\sigma_b=\sigma_m=1$, $t_{freq}=10$, and $s_{freq}=256$). The LSTM-Nudge is found to perform sufficiently well with this new observable quantity, as shown in Fig.~\ref{fig:au2} for the temporal modal coefficients and Fig.~\ref{fig:uu2} for the velocity field reconstruction. We also emphasize here that the similar behavior between observing $u(x,t)$ and observing $u^2(x,t)$ might be attributed to the values of the velocity field varying between $0$ and $1$. Thus, observing either the velocity or its squared value basically has the same pattern and range. For different situations, where the observable has significantly different pattern, the behavior might vary as well. For example, for a sine wave moving between $-1$ and $1$, observing the squared value (or the absolute value) would result in measurements between $0$ and $1$, neglecting the negative part.

\section{Concluding Remarks} \label{sec:conc}
In the current study, we have developed a methodology to utilize machine learning to cure model deficiency through online measurement data adopting ideas from dynamic data assimilation. In particular, an LSTM architecture has been trained to nudge prior predictions towards true state values using a combination of background information with sparse and noisy observations. The proposed framework is distinguished from previous studies in the sense that it is built on the assumption that all the computing ingredients are imperfect, including a truncated GROM model, erroneous initial conditions, and defective sensors. We have applied the proposed LSTM-Nudge to the 1D Burgers problem with a moving discontinuity, and investigated the effects of measurement noise and initial condition perturbation on its behavior. Although the framework works sufficiently well for a wide range of noise and perturbation, numerical experiments have indicated relatively more dependence of performance on measurement quality (noise). Meanwhile, we have found that sensors sparsity has minimal effects on results. We emphasize that the proposed framework represents one way of hybridizing human knowledge, physics-based models, measurement information, and data-driven tools to maximize their benefits rather than discarding any of them. This might represent a viable key enabler for the emerging digital twin applications. Nonetheless, the scalability of the approach has yet to be tested using more complex and higher-dimensional problems. 

\section*{Acknowledgments}
This material is based upon work supported by the U.S. Department of Energy, Office of Science, Office of Advanced Scientific Computing Research under Award Number DE-SC0019290. O.S. gratefully acknowledges their support. 
Disclaimer: This report was prepared as an account of work sponsored by an agency of the United States Government. Neither the United States Government nor any agency thereof, nor any of their employees, makes any warranty, express or implied, or assumes any legal liability or responsibility for the accuracy, completeness, or usefulness of any information, apparatus, product, or process disclosed, or represents that its use would not infringe privately owned rights. Reference herein to any specific commercial product, process, or service by trade name, trademark, manufacturer, or otherwise does not necessarily constitute or imply its endorsement, recommendation, or favoring by the United States Government or any agency thereof. The views and opinions of authors expressed herein do not necessarily state or reflect those of the United States Government or any agency thereof.
\bibliography{ref}

\end{document}